\documentclass{siamltex}
\usepackage{amsfonts}
\usepackage{amsmath}
\usepackage{amssymb}
\usepackage{array}
\usepackage{stmaryrd}

\def\dx{\hspace{2pt}{\rm d}x}
\def\d{\hspace{2pt}{\rm d}}

\def\l2{_{L_2(\Omega)}}
\def\nl2{_{[L_2(\Omega)]^n}}
\def\hmin{\underline{h}}

\def\E{\mathcal{E}}

\def\diam{{\rm{diam}}}

 \def\XXint#1#2#3{{\setbox0=\hbox{$#1{#2#3}{\int}$} 
\vcenter{\hbox{$#2#3$}}\kern-.5\wd0}}

\def\T{\mathcal{T}}

\def\R{\mathcal{R}}

\def\E{\mathcal{E}}

\newcommand{\sobz}[2]{\ensuremath{{\overset{\smash{\scriptscriptstyle\circ}}W}{}^{#1}_{#2}}}

\newcounter{saveeqn}

\newtheorem{prop}[theorem]{Proposition} 
\newtheorem{remark}[theorem]{\it Remark}

\numberwithin{equation}{section}

\title{A posteriori error estimates in the maximum norm for parabolic problems \thanks{The first author was partially supported by National Science Foundation grant DMS-0303378.  The second author's research has been partially supported by the Nuffield Foundation's "Young Researcher's Grant".}
}

\author{Alan Demlow\thanks{University of Kentucky, Department of Mathematics, 715 Patterson Office Tower, Lexington, Kentucky 40506--0027 ({\tt demlow@ms.uky.edu})} \and Omar Lakkis\thanks{Department of Mathematics, University of Sussex, Brighton, England, UK-BN1 9RF,  United Kingdom ({\tt o.lakkis@sussex.ac.uk})} \and  Charalambos  Makridakis\thanks{Department of Applied Mathematics,
University of Crete,
GR-71409 Heraklion, Greece
and
Institute for Applied and Computational Mathematics,
Foundation for Research and Technology-Hellas,
Vasilika Vouton P.O.Box 1527,
GR-71110 Heraklion, Greece
({\tt makr@tem.uoc.gr})}.}

\begin{document}

\maketitle

\renewcommand{\thefootnote}{\fnsymbol{footnote}}

\begin{abstract} We derive a posteriori error estimates in the $L_\infty((0,T];L_\infty(\Omega))$ norm for approximations of solutions to linear parabolic equations.  Using the elliptic reconstruction technique introduced by Makridakis and Nochetto and heat kernel estimates for linear parabolic problems, we first prove a posteriori bounds in the maximum norm for semidiscrete finite element approximations.  We then establish a posteriori bounds for a fully discrete backward Euler finite element approximation.  The elliptic reconstruction technique greatly simplifies our development by allowing the straightforward combination of heat kernel estimates with existing elliptic maximum norm error estimators.
 \end{abstract}


\begin{AM} 65N30
\end{AM}

\pagestyle{myheadings}
\thispagestyle{plain}
\markboth{A. DEMLOW, O. LAKKIS, AND C. MAKRIDAKIS}{MAXIMUM NORM ESTIMATES FOR PARABOLIC PROBLEMS}


\section{Introduction} \label{sec1}
We consider finite element approximations to the problem
\begin{equation}
\begin{split}
u_t-\Delta u&= f \hbox{ in } \Omega \times (0,T],
\\u & =  0 \hbox{ on } \partial \Omega \times [0,T],
\\u(x,0) & =  u_0(x).
\end{split} 
\label{eq1-1-1}
\end{equation}
Here $\Omega \subset \mathbb{R}^n$ ($n=2,3$) is a bounded polyhedral domain, $f$ is sufficiently smooth, $u_0 \in L_\infty(\Omega)$, and $u$ is a weak solution to (\ref{eq1-1-1}) lying in $L_\infty(0,T;L_\infty(\Omega))\cap H^1(0,T;H^{-1}(\Omega))$ (a subset of $C^0(0,T;L_\infty)$.  

Adaptive finite element methods for approximating solutions to parabolic partial differential equations are popular because of their ability to efficiently resolve singularities and other rapid local variations in solutions.  While most adaptive finite element methods are designed to control only energy norms of solutions, in many applied problems the goal output of a finite element computation is related to some other norm or functional of the solution.  In this work we address control of the {\it maximum} error $\|u-u_h\|_{L_\infty(\Omega \times [0,T])}$ for finite element approximations $u_h$ of $u$.  Ensuring good pointwise approximation of $u$ is natural in many situations where $u$ represents some physical quantity.  Pointwise error control is also a natural goal when computing free boundaries, for example via level set methods (cf. \cite{DDE05}).  Several recent papers have addressed adaptive finite element methods for controlling pointwise errors in elliptic problems (cf. \cite{Noc95}, \cite{DDP00}, \cite{NSV03}, \cite{NSV05}, \cite{NSSV06}, \cite{De06}, \cite{Dem07}).  However, the only previous pointwise a posteriori estimates for parabolic problems that we are aware of are contained in \cite{EJ95} and \cite{Bo00}, which we describe below.     

The goal of this work is to prove a posteriori error estimates in $L_\infty(\Omega \times [0,T])$ for semi- and fully-discrete finite element approximations to (\ref{eq1-1-1}).  For practical purposes, our main result is an easily-computable error estimator for the backward Euler finite element discretization of (\ref{eq1-1-1}).  In order to describe this estimate, we introduce some definitions and notation.  Let $0=t_0<t_1<...<t_N=T$, $I_i=(t_{i-1}, t_i)$, and $\tau_i=t_i-t_{i-1}$.  For each $0 \le i \le N$, let $\T_i$ be a triangulation of $\Omega$.  We place only standard restrictions on the triangulations, requiring in particular that all triangles have aspect ratios that are uniformly bounded with respect to $ i=0, \dotsc, N$ and that the triangulations are ``edge-to-edge'' (i.e., hanging nodes are not allowed).  Let $S_0^i$ be a finite element space consisting of the functions that are continuous piecewise polynomials of degree $k$ on $\T_i$ and which are $0$ on $\partial \Omega$.  Letting $v^i(x)=v(t_i,x)$ for any function $v$ defined on $\Omega \times [0,T]$, we discretize the weak form of (\ref{eq1-1-1}) by letting $u_h^0 \in S_0^0$ approximate $u_0$ and then defining $u_h^i \in S_0^i$, $1 \le i \le N$, via the implicit Euler recursion
\begin{equation}
\frac{1}{\tau_i} \int_\Omega (u_h^i-u_h^{i-1}) \phi_i \dx+ \int_\Omega \nabla u_h^i \nabla \phi_i \dx=\int_\Omega f^i \phi_i \dx ~ \hbox{ for all } \phi_i \in S_0^i.
\label{eq1-1-2}
\end{equation}
In addition, let $g^i=f^i-\frac{u_h^i-u_h^{i-1}}{\tau_i}$, $i \ge 1$.  The definition of $g^0$ is slightly different; cf. \S 4.2.  Then for any $1 \le j \le N$,
\begin{equation}
\begin{split}
\|u^j-u_h^j\|_{L_\infty(\Omega)} \le & \|u_0-u_h^0\|_{L_\infty(\Omega)}
\\ & +C(\Omega) (2+ c(n) \ln \frac{t_j}{\tau_j}) (\ln \hmin)^2 \max_{0 \le i \le j} \E_{\infty,0}(u_h^i, g^i)
\\ &+ \sum_{i=1}^j \int_{I_i} \|f-f^{i}\|_{L_\infty(\Omega)} \d t+\frac{\tau_i}{2}\|g^i-g^{i-1}\|_{L_\infty(\Omega)}.
\end{split}
\label{eq1-1-3}
\end{equation}
Here $c(n)=\frac{3^n}{2^{n/2+1}}$, $\E_{\infty,0}(u_h^i,g^i)$ is a standard and easily-computable residual-type estimator depending only on $u_h$, $g^i$, and the mesh $\T_i$, and $\hmin$ is the minimum diameter of elements lying in $\cup_{i=1}^N \T_i$.  A more precise definition of $\E_{\infty,0}$ is provided in \S \ref{sec2-2},  and Theorem \ref{t4-4} contains a precise statement of results.  

The error bound on the right hand side of (\ref{eq1-1-3}) consists of:
\begin{itemize}
\item an {\it initial data estimator} $\|u_0-u_h^0\|_{L_\infty(\Omega)}$;
\item a {\it spatial estimator} $C(\Omega) (2+c(n) \ln \frac{t_j}{\tau_j}) (\ln \hmin)^2 \max_{0 \le i \le j} \E_{\infty,0}(u_h^i, g^i)$ which accounts for spatial errors; and
\item  a {\it time estimator} $\sum_{i=1}^j \int_{I_i} \|f-f^{i}\|_{L_\infty(\Omega)} \d t+\frac{\tau_i}{2}\|g^i-g^{i-1}\|_{L_\infty(\Omega)}$ which accounts temporal errors.
\end{itemize}
Note that the constants in front of the initial data and time estimators are 1, so that estimating these terms does not require estimating unknown constants.  While the spatial estimator does contain the constant $C(\Omega)$, we shall show later that $C(\Omega)$ depends only on properties of an underlying {\it elliptic} a posteriori error estimator.  

In addition to the bound (\ref{eq1-1-3}), we also establish several other a posteriori estimates for semidiscrete finite element approximations of $u$ as well as for the fully discrete scheme (\ref{eq1-1-2}).  While more difficult than (\ref{eq1-1-3}) to employ in practical adaptive codes, these estimates provide additional insight into a posteriori theory for pointwise norms, for example by establishing that the spatial terms in the middle line of (\ref{eq1-1-3}) can sometimes be bounded instead in a weaker negative norm.

We next briefly compare the estimate (\ref{eq1-1-3}) with the results of the previously cited works \cite{EJ95} and \cite{Bo00}.   \cite{EJ95} contains pointwise a posteriori estimates for parabolic problems discretized in space using standard piecewise linear finite element schemes and in time using a discontinuous Galerkin approach.    However, these estimates are stated without proof, do  not appear to be proven elsewhere in the literature, and additionally assume restrictive hypotheses on the spatial finite element meshes.    \cite{Bo00} similarly employs a discontinuous Galerkin time discretization and piecewise linear finite element spatial discretization.   The proofs in \cite{Bo00} involve proving quasi-optimal regularity estimates for a regularization of the parabolic Green's function, which is fairly difficult and also leads to an uncomputable a priori term in the upper bound.  This method of proof essentially involves imitating in a parabolic context the maximum norm estimates for elliptic problems originally proven in \cite{Noc95} and \cite{DDP00}.  We finally note that the results of both \cite{Bo00} and \cite{EJ95} are restricted to convex polyhedral domains.  

We emphasize some features of (\ref{eq1-1-3}) that contrast positively with existing results.  First, in the present work $\Omega$ may be a nonconvex polyhedral domain (including a domain with cracks).  In addition, we allow arbitrary orders of finite element spaces, and in particular do not restrict ourselves to piecewise linear elements.  The estimate (\ref{eq1-1-3}) also does not require any impractical restrictions on the spatial mesh (in contrast to \cite{EJ95}), and does not contain any uncomputable terms depending on $u$ in the upper bound (in contrast to \cite{Bo00}).  Finally, as we discuss further below, the proof of (\ref{eq1-1-3}) is quite straightforward because we are able to {\it reuse} difficult elliptic results instead of {\it imitating} their proofs in a parabolic context.  Thus the results that we present here are to our knowledge the first rigorously proven, fully a posteriori pointwise error estimates for finite element methods for parabolic problems.

Essential to our development is the {\it elliptic reconstruction} technique introduced in \cite{MN03} in the context of semidiscrete problems and extended to fully discrete problems in \cite{LM06}.  In essence, the elliptic construction $\R u_h$ is a continuous elliptic representation of the discrete solution $u_h$, and   $u_h$ is the elliptic finite element approximation to $\R u_h$ with respect to the finite element space under consideration.  Thus any a posteriori error estimates which are valid for elliptic problems on $\Omega$ may be used to estimate $\R u_h-u_h$.  The overall error $u-u_h$ may then be bounded by first estimating $u-\R u_h$ using PDE techniques for continuous parabolic problems and then estimating $\R u_h-u_h$ using elliptic a posteriori estimators.  The elliptic reconstruction may thus be regarded as an a posteriori analogue to the Ritz-Wheeler projection in standard a priori error analysis for parabolic problems (cf. \cite{Wh73}, \cite{Th97}).  

Our use of the elliptic reconstruction technique in the context of pointwise error estimation for parabolic problems highlights its ability to fully leverage existing elliptic estimates.  In particular, establishing a rigorous theory for a posteriori estimation of pointwise errors for Poisson's problem on polyhedral domains was a technically difficult enterprise (cf \cite{Noc95}, \cite{DDP00}, \cite{NSSV06}).  Relying on these elliptic results instead of mimicking them, our proofs employ only basic estimates for the heat kernel for the {\it continuous} problem (\ref{eq1-1-1}) and are relatively short and straightforward.  A direct ``parabolic'' approach to the problem which does not use the elliptic reconstruction is also possible (cf \cite{Bo00}), but such an approach is much more technically involved and as already mentioned has not led to the sharpest possible results.

An outline of the paper follows.  In \S2, we provide common preliminaries and recall some facts concerning residual-type a posteriori error estimation for elliptic problems.  In \S3 we prove a posteriori estimates for semidiscrete approximations of (\ref{eq1-1-1}), while in \S4 we consider a backward Euler time discretization of (\ref{eq1-1-1}).  

\section{Preliminaries}
In this section we provide a number of preliminaries regarding heat kernel estimates, maximum norm a posteriori estimates for elliptic problems, and issues concerning mesh compatibility that arise in some of our estimates for fully discrete schemes.  

\subsection{Notation}
We begin by defining suitable notation. $W_p^j(\Omega)$, where $j$ is a nonnegative integer,  will denote the standard Sobolev space of functions having $j$ derivatives in $L_p(\Omega)$, and $\sobz{1}{p}(\Omega)$ will denote the functions in $W_p^1(\Omega)$ which in addition have zero trace on the boundary $\partial \Omega$.  In addition, $(\cdot, \cdot)$ denotes either the $L_2(\Omega)$- or $[L_2(\Omega)]^n$-inner product.  Finally, we denote by $L_p([a,b],W_q^j(\Omega))$ the functions whose spatial $W_q^j$ norm lies in $L_p$ over the time interval $[a,b]$.  

\subsection{Analytical preliminaries}
Our analysis relies heavily on properties of the heat kernel for the problem (\ref{eq1-1-1}).  We sum up the necessary results in the following lemma.
\begin{lemma}
\label{lem2-1}
Let $\Omega \subset \mathbb{R}^n$ be a bounded open domain.  Then there exists a Green's function $G(x,t; y, s)$ for the problem (\ref{eq1-1-1}).  That is, there exists a kernel $G$ such that for $u$ satisfying (\ref{eq1-1-1}) and $(x,t) \in \Omega \times (0,T]$, 
\begin{equation}
u(x,t)=\int_\Omega G(x,t; y,0) u_0(y)  \d y+\int_0^t \int_\Omega G(x,t; y,s) f(y,s) \dx \d t
\label{eq2-3-1}
\end{equation}
 is a weak solution of (\ref{eq1-1-1}).  For $s<t$, $G$ additionally satisfies the bound
 \begin{equation}
\|G(x,t; \cdot,s)\|_{L_1(\Omega)} \le 1.
\label{eq2-3-4}
\end{equation}
 
  
  Let also $2 <p,q  \le \infty$ with 
  \begin{equation}
  \frac{n}{2p}+\frac{1}{q}<\frac{1}{2}.  
\nonumber
  \end{equation}
  Then we have in addition that $G(x,t; \cdot, \cdot) \in L_{q'}([0,T], \sobz{1}{p'}(\Omega))$, where $p'$ and $q'$ are the conjugate exponents to $p$ and $q$.  Also, 
  \begin{equation}
\|G(x,t; \cdot, \cdot)\|_{L_{q'}([0,T], W_{p'}^1(\Omega))} \le C_{p,q}(T),
\label{eq2-3-3}
\end{equation}
where $C_{p,q}$ depends on $p$, $q$, $|\Omega|$, and $T$.  In addition, $G(x,t;\cdot,s) \in H_0^1(\Omega)$ for $0 \le s<t$.

Finally, there exists a constant $c(n)$ depending only on the space dimension $n$ such that for $s<t$, \begin{equation}
\|G_s(x,t; \cdot, s)\|_{L_1(\Omega)} \le \frac{c(n)}{t-s}.
\label{eq2-3-5}
\end{equation}
Here we use the notation $G_s(x,t; \cdot, s)=\frac{\partial}{\partial s} G(x,t; \cdot, s)$.  
\end{lemma}

{\it Proof: }  
The existence of a Green's function satisfying (\ref{eq2-3-1}) and (\ref{eq2-3-3}) may be found in Theorem 6 (p. 657) and Theorem 9 (p. 671) of the fundamental work \cite{Ar68} of Aronson.  To prove (\ref{eq2-3-4}), we note that 
\begin{equation}
0 \le G(x,t;\cdot ,s) \le \Gamma(x,t; y,s)=(4 \pi (t-s))^{-n/2} e^{-\frac{|x-y|^2}{4(t-s)}}.
\label{eq2-3-6}
\end{equation}
That is, the heat kernel on $\Omega$ is bounded pointwise by the heat kernel on $\mathbb{R}^n$.  This fact may be proven e.g. by combining Lemma 7 (p. 677) of \cite{Ar68} with Corollary 8.2 and Theorem 8.3 of \cite{Dan00}.  Inequality (\ref{eq2-3-4}) then follows from the fact that $\int_{\mathbb{R}^n} \Gamma(x,t; y,s) \d y=1$ for $t<s$.  

In order to prove (\ref{eq2-3-5}), we apply Corollary 5 of \cite{Dav97} with $\delta=\frac{1}{2}$ and $\epsilon=\frac{1}{9}$ to find that 
\begin{equation}
G_s(t,x; y,s) \le 2^{n/2-1} (t-s)^{-1} (4 \pi (t-s))^{-n/2} e^{-\frac{|x-y|^2}{9(t-s)}}.
\label{eq2-3-7}
\end{equation}
Performing a change of variables and integrating over $\mathbb{R}^n$ yields (\ref{eq2-3-5}) with $c(n)=\frac{3^n}{2^{n/2+1}}$.  
\hfill $\Box$

\begin{remark}{\rm
One may take advantage of the bound (\ref{eq2-3-6}) and (\ref{eq2-3-7}) in order to explicitly incorporate dissipation of the heat kernel into (\ref{eq2-3-4}) and (\ref{eq2-3-7}).  For example, for $n=2$ and $s<t$ one may compute that
\begin{eqnarray}
\|G(x,t; \cdot, s)\|_{L_1(\Omega)} &\le& 1-e^{-\frac{\diam(\Omega)^2}{4(t-s)}},
\label{eq2-3-8}
\\ \|G_s(x,t:,\cdot, s)\|_{L_1(\Omega)} &\le& \frac{c(n)}{t-s} (1-e^{-\frac{\diam(\Omega)^2}{9(t-s)}}).
\label{eq2-3-9}
\end{eqnarray}
It is possible to incorporate (\ref{eq2-3-8}) and (\ref{eq2-3-9}) into pointwise a posteriori bounds with no great difficulty, and we shall pursue this briefly in Corollary \ref{cor4-4}.  
} \end{remark}

\subsection{Elliptic a posteriori estimates} \label{sec2-2}
In this section we cite several results that will enable us to bound a posteriori the elliptic reconstruction terms appearing in our estimates for parabolic problems.  In this subsection we assume that $v$ satisfies
\begin{equation}
\begin{split}
-\Delta v &=g \hbox{ in } \Omega,
\\ v & =  0 \hbox{ on } \partial \Omega
\end{split}
\nonumber
\end{equation}
where $\Omega\subset \mathbb{R}^n$, $n=2,3$ is a polyhedral domain.  We additionally assume that $\T$ is a shape-regular simplicial decomposition of $\Omega$, and define the Lagrange finite element space $S=\{ w_h \in H^1(\Omega): w_h |_{K} \text{ is a polynomial of degree } k \text{ on } K, ~K \in \T \}$.  Let also  $S_0=S \cap H_0^1(\Omega)$.  Finally, let $v_h \in S_0$ be the finite element approximation to $v$ defined by
\begin{equation}
\int_\Omega \nabla v_h \nabla w_h \dx = \int_\Omega f w_h \dx ~ \hbox{ for all } w_h \in S_0.
\nonumber
\end{equation}

Our parabolic results assume a posteriori bounds for $\|v-v_h\|_{L_\infty(\Omega)}$, and in some circumstances also for $\|v-v_h\|_{W_p^{-1}(\Omega)}$.  With $\frac{1}{p}+\frac{1}{p'}=1$, here 
\begin{equation}
\|w\|_{W_p^{-1}(\Omega)}=\sup\{(w,x):{z \in \sobz{1}{p'}(\Omega), \|z\|_{W_{p'}^1(\Omega)}=1} \}.
\nonumber
\end{equation}
We shall employ residual-type estimates.  We first define the jump residual $\llbracket \nabla v_h \rrbracket$ on an $(n-1)$-dimensional element face $e=K_1 \cap K_2$, where $K_1, K_2 \in \T$.  Let $\vec{n}$ be an arbitrary unit normal vector on $e$, and for $x \in e$ let 
\begin{equation}
\llbracket \nabla v_h \rrbracket(x)=\lim_{\delta \to 0}  (\nabla v_h(x+\delta \vec{n})-\nabla v_h(x-\delta \vec{n})) \cdot \vec{n}.   
\nonumber
\end{equation}
Let also $h_K$ be the diameter of the element $K$.  For $1 \le p \le \infty$ and $j \ge 0$, we then define the elementwise error indicator
\begin{equation}
\eta_{p,-j}(K)=h_K^{2+j} \|g+\Delta v_h\|_{L_p(K)}+h_K^{j+1+\frac{1}{p}}\|\llbracket \nabla v_h \rrbracket \|_{L_p(\partial K)}.
\nonumber
\end{equation}
Finally, we define the global estimator
\begin{equation}
\E_{p,-j}(v_h, g)=\left \{ \begin{array}{l} (\sum_{K \in \T} \eta_{p,-j}(K)^p)^{1/p}, ~1 \le p < \infty, \\ \max_{K \in \T} \eta_{\infty,-j}(K), ~p=\infty. \end{array} \right .
\label{eq2-4-4}
\end{equation}

We first quote an elliptic pointwise estimate which holds for all orders of finite element spaces and all polyhedral domains; cf. \cite{NSSV06} for a proof.  
\begin{lemma}
\label{lem2-4-1}
Assume that $\Omega$ is an arbitrary polyhedral domain in $\mathbb{R}^n$, $n=2,3$, and let $\hmin=\min_{K \in \T} h_K$.  Then
\begin{equation}
\|v-v_h\|_{L_\infty(\Omega)} \le C(\Omega) (\ln \hmin)^2 \E_{\infty,0}(v_h, g).
\nonumber
\end{equation}
\end{lemma}

When allowed by the domain geometry, it may be advantageous to instead measure the error in a negative norm.  In such cases we shall employ the following result.
\begin{lemma}
\label{lem2-4-2}
Assume that $\Omega$ is a convex polyhedral domain in $\mathbb{R}^n$, $n=2,3$, where the maximum  vertex opening angle (for $n=2$) or edge opening angle (for $n=3$) is denoted by $\omega=\frac{\pi}{\alpha}$, $\alpha>1$.  Assume also that the degree of the finite element space $S_0$ is at least 2, that is, $S_0$ contains the continuous piecewise quadratic functions.  Then for $\frac{2}{\alpha-1} < p < \infty$,
\begin{equation}
\|v-v_h\|_{W_p^{-1}(\Omega)} \le C(p,\Omega) \E_{p,-1}(v_h,g).
\label{eq2-4-6}
\end{equation}
\end{lemma}
{\it Proof:}  
Assume that $-\Delta w=z$, where $z \in W_{p'}^1(\Omega)$ for $\frac{1}{p}+\frac{1}{p'}=1$.  Combining the comments of \S4.c and Corollary 3.9 of \cite{Da92} yields the regularity result
\begin{equation}
\|w\|_{W_{p'}^3(\Omega)} \le C(p',\Omega)\|z\|_{W_{p'}^1(\Omega)}.
\nonumber
\end{equation}
Given this estimate, the result (\ref{eq2-4-6}) may be obtained using a duality argument and standard techniques for proving residual-type a posteriori bounds. 
\hfill $\Box$

\begin{remark}{\rm We emphasize that Lemma \ref{lem2-4-2} only holds on {\it convex} polyhedral domains.  It should be possible to similarly prove usable negative norm estimates on nonconvex polygonal domains in $\mathbb{R}^2$, but explicit information about corner singularities would appear in such estimators (cf. \cite{LN03} for analogous global $L_2$ bounds on nonconvex polygonal domains).  On nonconvex polyhedral domains in $\mathbb{R}^3$, such a result is much less practical as the precise nature of vertex singularities is often difficult to ascertain. }
\end{remark}

\begin{remark}{\rm Care must be taken when employing (\ref{eq2-4-6}) as the constant $C(p,\Omega)$ degenerates as $p \rightarrow \infty$, and possibly as $p$ approaches the lower bound $\frac{2}{\alpha-1}$ as well.  In particular, one is not able to choose $p(\alpha)$ so that $C(p(\alpha),\Omega)$ remains bounded as the maximum edge opening angle approaches $\pi$ (i.e., as $\alpha \rightarrow 1$).  
}
\end{remark}

\subsection{Compatible meshes and estimates for differences in finite element solutions} \label{sec2-1a}
Some of our fully discrete a posteriori estimates require bounding elliptic finite element errors of the form $v_1-v_2-(v_{h1}-v_{h2})$, where $v_{h1}$ and $v_{h2}$ lie in different finite element spaces (in particular, in finite element spaces defined on meshes at adjacent time steps).  In this subsection we make  assumptions and definitions on the pair of meshes that allow us to establish such estimates.  We follow closely Appendices A and B of \cite{LM06}, so we shall be brief and refer the reader to that work for more details.  

Two simplicial decompositions  $\T_1$ and $\T_2$ of $\Omega$ are said to be {\it compatible} if they are derived from the same macro triangulation $\mathcal{M}$ by an {\it admissible} refinement procedure which preserves shape regularity and assures that for any elements $K \in \T_1$ and $K' \in \T_2$, either $K \cap K'=\emptyset$, $K \subset K'$, or $K' \subset K$.  The bisection-based refinement procedure used for example in the ALBERTA finite element toolbox (cf. \cite{SS05}) is known to be admissible.  

There is a natural partial ordering of compatible triangulations, with $\T_1 \le \T_2$ if $\T_2$ is a refinement of $\T_1$.  The {\it finest common coarsening} $\T_1 \wedge \T_2$ of $\T_1$ and $\T_2$ is defined in a natural way, and $\hat{h}=\max ( h_1, h_2)$, where $h_1$, $h_2$, and $\hat{h}$ are the local mesh size functions for $\T_1$, $\T_2$, and $\T_1 \wedge \T_2$.  Finally, let $S_1$ and $S_2$ be  finite element spaces of degree $k$ on $\T_1$ and $\T_2$.  $\hat{S}=S_1 \cap S_2$ is then the corresponding space of degree $k$ on $\T_1 \wedge \T_2$.  Standard interpolation inequalities hold for all of the above-mentioned spaces, though the constants in these bounds may depend on the number of refinement steps used to pass from $\T_1$ to $\T_2$.  We again refer to Appendix B of \cite{LM06}, Appendix B for more discussion.

Let $v_{h1} \in S_1$ and $v_{h2}\in S_2$ be the elliptic finite element approximations to $v_1$ and $v_2$ lying in finite element spaces $S_1$ and $S_2$ defined on different meshes $\T_1$ and $\T_2$, respectively.  Here we assume that $-\Delta v_1=g_1$, $-\Delta v_2=g_2$, and $v_1=v_2=0$ on $\partial \Omega$.   In essence, Lemma 2.3 and Lemma 2.4 still hold, but with the local mesh size $h$ replaced by the local mesh size $\hat{h}$ of the finest common coarsening.  Let $\Sigma_i$, $i=1,2$, be the union of the faces of elements lying in $\T_i$.  For $\hat{K} \in \T_1 \wedge \T_2$, let $\Sigma_{\hat{K}}=(\Sigma_1 \cup \Sigma_2) \cap \hat{K}$, where $\hat{K}$ is taken to be closed. 

For $1 \le p \le \infty$ and $j \ge 0$, we then define the elementwise error indicator
\begin{equation}
\hat{\eta}_{p,-j}(\hat{K})=\hat{h}_{\hat{K}}^{2+j} \|g_1-g_2+\Delta v_{h1}-\Delta v_{h2}\|_{L_p(\hat{K})}+h_{\hat{K}}^{j+1+\frac{1}{p}}\|\llbracket \nabla (v_{h1}-v_{h2}) \rrbracket \|_{L_p(\Sigma_{\hat{K}})}.
\nonumber
\end{equation}
Finally, we define the global estimator
\begin{equation}
\hat{\E}_{p,-j}(v_{h1}-v_{h2}, g_1-g_2; \T_1, \T_2)=\left \{ \begin{array}{l} (\sum_{\hat{K} \in \T_1 \wedge \T_2} \hat{\eta}_{p,-j}(\hat{K})^p)^{1/p}, ~1 \le p < \infty, \\ \max_{\hat{K} \in \T_1 \wedge \T_2} \hat{\eta}_{\infty,-j}(\hat{K}), ~p=\infty. \end{array} \right .
\nonumber
\end{equation}

\begin{corollary}  
\label{cor2-7}
Assume that $\Omega$ is an arbitrary polyhedral domain in $\mathbb{R}^n$, $n=2,3$, and that $\T_1$ and $\T_2$ are compatible triangulations.  In addition, define $\hat{\hmin}=\min_{x \in \Omega} \min(h_1(x), h_2(x))$.  Then
\begin{equation}
\|v_1-v_2-(v_{h1}-v_{h2})\|_{L_\infty(\Omega)} \le C(\Omega) (\ln \hat{\hmin})^2 \hat{\E}_{\infty,0}(v_{h1}-v_{h2}, g_1-g_2; \T_1, \T_2).
\label{eq2-4-10}
\end{equation}
Here $C(\Omega)$ depends on the number of refinement steps used to pass from $\T_1$ to $\T_2$.
\end{corollary}

\begin{corollary}
\label{cor2-8}
Assume that $\Omega$ is a convex polyhedral domain in $\mathbb{R}^n$, $n=2,3$, where the maximum  vertex opening angle (for $n=2$) or edge opening angle (for $n=3$) is denoted by $\omega=\frac{\pi}{\alpha}$, $\alpha>1$.  Assume also that the degree $k$ of the finite element spaces $S_1$ and $S_2$ is at least 2, that is, both spaces contain the continuous piecewise quadratic functions. Assume also that $\T_1$ and $\T_2$ are compatible triangulations. Then for $\frac{2}{\alpha-1} < p < \infty$,
\begin{equation}
\|v_1-v_2-(v_{h1}-v_{h2})\|_{W_p^{-1}(\Omega)} \le C(p,\Omega) \hat{\E}_{p,-1}(v_{h1}-v_{h2},g_1-g_2; \T_1, \T_2).
\label{eq2-4-11}
\end{equation}
Here $C(p,\Omega)$ depends on $p$ and the number of refinement steps used to pass from $\T_1$ to $\T_2$.
\end{corollary}

The main observation used to derive the above corollaries is the fact that $(v_1-v_2)-(v_{h1}-v_{h2})$ is Galerkin orthogonal to the space $S_1 \cap S_2$.  The compatibility of $\T_1$ and $\T_2$ ensures that this intersection is rich enough to obtain (\ref{eq2-4-10}) and (\ref{eq2-4-11}).  The proofs otherwise follow closely those of Lemma \ref{lem2-4-1} and Lemma \ref{lem2-4-2}, and we do not give details here.  

\section{Analysis of the semidiscrete scheme}

\subsection{Semidiscrete finite element approximation}

For simplicity in handling finite element approximations, we assume that $\Omega$ is a polyhedral domain in $\mathbb{R}^n$, $n=2,3$, and that $\T_h$ is a simplicial decomposition of $\Omega$.  We emphasize that we admit here general polyhedral domains, including those having cuts or cracks (as described in the standard reference \cite{Da88}).  Let $S \subset H^1(\Omega)$ and $S_0 =S\cap H_0^1(\Omega)$ be standard simplicial Lagrange finite element spaces as in \S \ref{sec2-2}.  The semidiscrete approximation $u_h \in C([0,T], S_0)$ of $u$ then satisfies
\begin{equation}
\begin{split}
(u_{h,t},v_h)+(\nabla u_h, \nabla v_h)=&(f,v_h), ~ v_h\in S_0 ~ and~t \in (0,T],
\\ u_h(0)=& u_h^0,
\end{split}
\label{eq3-1-1}
\end{equation}
where $u_h^0 \in S_0$ is a finite element approximation to $u_0$.    We let $P_h:L_2 \rightarrow S$ be the $L_2$ projection onto the finite element space $S$, and additionally define the modified discrete Laplacian $-\Delta_h(t) :H_0^1(\Omega) \rightarrow (S_0+P_h f(t))$ by
\begin{equation}
(\nabla u, \nabla v_h)=(-\Delta_h(t) u, v_h), ~ v_h \in S_0.
\label{eq3-1-1a}
\end{equation}
From (\ref{eq3-1-1}) and (\ref{eq3-1-1a}), we have the pointwise formula
\begin{equation}
-\Delta_h u_h=P_h f-u_{h,t}.
\label{eq3-1-1b}
\end{equation}

\begin{remark}{\rm
Our definition of $-\Delta_h$ is nonstandard in that here $-\Delta_h u$ has nonzero boundary values which depend on the data $f$ in addition to $u$.  We use this definition in order to maintain consistency in the pointwise relationship (\ref{eq3-1-1b}).  In particular, because $u_{h,t}=0$ on $ \partial \Omega$, we also have $P_h f=-\Delta u_h$ on $\partial \Omega$.  Note that we instead could enforce this relationship by letting $-\Delta_h u_h \in S_0$ and taking the $L_2$ projection of $f$ onto $S_0$.  This distinction will make little practical difference in our development, but it is possible to define the elliptic reconstruction in such a way that $-\Delta_h u_h$ and $P_h f$ must be computed a posteriori (cf. \cite{LM06}).   }\end{remark}

\subsection{Elliptic reconstruction for the semidiscrete problem}

Given a finite element approximation $u_h$, we define its elliptic reconstruction $\R u_h \in H_0^1(\Omega)$ by 
\begin{equation}(\nabla \R u_h, \nabla v)=(g,v), ~ v \in H_0^1(\Omega),
\label{eq3-1-2}
\end{equation}
where $g=f-u_{h,t}$.  We thus have $-\Delta _h u_h=P_h f-u_{h,t}$, $-\Delta \R u_h=f-u_{h,t}$, and $-\Delta u=f-u_t$.  Note that $\R=\R(t)$ is a time-dependent operator, but we shall depress its dependence on $t$ in the sequel in order to avoid unnecessary clutter in our notation.  We will deal more explicitly with the time dependence of the reconstruction operator in our analysis of the fully discrete scheme.  Also note that we may differentiate (\ref{eq3-1-2}) with respect to $t$ to obtain
\begin{equation}
(\nabla (\R u_h)_t, \nabla v)=(g_t,v), ~ v \in H_0^1(\Omega),
\nonumber
\end{equation}
where $g_t=f_t-u_{h,tt}$.

Combining (\ref{eq3-1-1a}), (\ref{eq3-1-1b}), and (\ref{eq3-1-2}), we find that $u_h$ and $\R u_h$ satisfy the Galerkin orthogonality relationship 
\begin{equation}
(\nabla(\R u_h-u_h),\nabla  v_h)=0, ~v_h \in S_0.
\label{eq3-1-3}
\end{equation}
Differentiating (\ref{eq3-1-3}) with respect to $t$ also yields the time-differentiated Galerkin orthogonality relationship
\begin{equation}
(\nabla (\R u_h-u_h)_t, \nabla v_h)=0, ~ v_h \in S_0(\Omega).
\nonumber
\end{equation}

In addition, it is easy to calculate that for $0 < t \le T$ and $v \in H_0^1(\Omega)$,
\begin{equation} 
((u-\R u_h)_t, v)+(\nabla(u-\R u_h),\nabla v)=((u_h -\R u_h)_t, v).
\label{eq3-1-3a}
\end{equation}
The fact that $u -\R u_h$ thus satisfies a parabolic equation will play a fundamental role in our development.  

\begin{remark}{\rm 
In \cite{MN03}, the elliptic reconstruction is defined by 
\begin{equation}
(\nabla \R u_h, \nabla v)=(-\Delta_h u_h+f-P_h f, v), ~ v \in H_0^1(\Omega).
\label{eq3-1-4}
\end{equation}
From (\ref{eq3-1-1b}), we see that $-\Delta_h u_h+f-P_h f=f-u_{h,t}$ so that (\ref{eq3-1-4}) and (\ref{eq3-1-2}) are equivalent.  In fact, the elliptic reconstruction allows us to write the pointwise form (\ref{eq3-1-1b}) of the discrete equation as 
\begin{equation}
u_{h,t}-\Delta \R u_h =f.
\nonumber
\end{equation}
The above equation does not involve the discrete Laplacian and thus allows for a straightforward comparison with the PDE (\ref{eq1-1-1}), leading in its weak form to (\ref{eq3-1-3a}).  We use the definition (\ref{eq3-1-2}) because in what follows we employ residual estimators to estimate the elliptic error $\R u_h -u_h$.  These estimators  require pointwise access to the right-hand-side data for $\R u_h$.  It is not practical to directly compute $-\Delta_h u_h$ or $P_h f$, and in \cite{MN03}, the authors develop an expression for the residual that does not involve the operators $-\Delta_h$ or $P_h$ (cf. p. 1592).  Thus the definition in \cite{MN03} emphasizes the underlying structure of the reconstruction operator, that is, $\R=(-\Delta)^{-1} (-\Delta_h)$ up to terms that are $L_2$-orthogonal to the finite element space.   Our equivalent definition instead reflects the practical concern of computing using the resulting a posteriori error estimates.  


} \end{remark}

\subsection{Semidiscrete reconstruction results}

Our a posteriori estimates are based upon the following theorem.
\begin{theorem}
\label{t3-1}
Let the assumptions and definitions of \S3.1 and \S3.2 hold.  Then for any $0 < t_0 \le T$, 
\begin{equation}
\begin{split}
\|(u-u_h)(t_0)\|_{L_\infty(\Omega)} \le& \|(\R u_h-u_h)(t_0)\|_{L_\infty(\Omega)}+\|(u-\R u_h)(0)\|_{L_\infty(\Omega)} 
\\&+\|(u_h -\R u_h)_t \|_{L_1([0,t_0]; L_\infty(\Omega))}. 
\end{split}
\label{eq3-2-2}
\end{equation}
Alternatively, let $2 < p,q \le \infty$ satisfy $\frac{n}{2p}+\frac{1}{q}<\frac{1}{2}$.  Then 
\begin{equation}
\begin{split}
\|(u-u_h)(t_0)\|_{L_\infty(\Omega)} \le& \|(\R u_h-u_h)(t_0)\|_{L_\infty(\Omega)}+\|(u-\R u_h)(0)\|_{L_\infty(\Omega)}
\\&+C_{p,q}(t_0) \|(u_h -\R u_h)_t\|_{L_q([0,t_0]; W_p^{-1}(\Omega))}.
\end{split}
\label{eq3-2-1}
\end{equation}
\end{theorem}
{\it Proof:}
For any $x_0 \in \Omega$, 
\begin{equation}
|(u-u_h)(x_0,t_0)| \le |(u-\R u_h )(x_0,t_0)|+\|(\R u_h -u_h)(t_0)\|_{L_\infty(\Omega)}.  
\nonumber
\end{equation}
Using (\ref{eq2-3-1}) and (\ref{eq3-1-3a}), we find that
\begin{equation}
\begin{split}
(u-\R u_h)(x_0,t_0)=&\int_\Omega G(x_0,t_0; y,0) (u-\R u_h)(y,0) \d y
\\ & +\int_0^{t_0} \int_\Omega G(x_0, t_0; y,s) (u_h -\R u_h)_t (y,s) \d s.
\end{split}
\nonumber
\end{equation}
The first term on the right hand side above may be bounded by combining H\"older's inequality with (\ref{eq2-3-4}).  In order to bound the second term, we apply (\ref{eq2-3-4}) to obtain (\ref{eq3-2-2}) or  (\ref{eq2-3-3}) to obtain (\ref{eq3-2-1}).  \hfill $\Box$


\begin{remark}{\rm  We have assumed a polyhedral domain and a specific type of finite element space in Theorem \ref{t3-1}, but similar results hold under more general circumstances.  Indeed, the analytical results of Lemma \ref{lem2-1} hold on general bounded domains, and only those estimates along with the relationship (\ref{eq3-1-3a}) are used in the proof.  The bound (\ref{eq2-3-4}) also holds for a fairly general class of elliptic differential operators, though in more general cases one must perhaps replace 1 on the right hand side by an unknown constant with unknown dependence on $T$.  The reconstruction technique is thus able to transfer most  difficulties and issues concerning the precise nature of the finite element approximation (including for example the type of elements used and  difficulties arising from finite element approximations on nonpolygonal domains) to the a posteriori estimation of elliptic errors.  
}
\end{remark}

\subsection{A posteriori error estimates for the semidiscrete problem}
In this subsection we shall estimate the right hand sides of  (\ref{eq3-2-2}) and (\ref{eq3-2-1}) a posteriori using the residual estimators of \S\ref{sec2-2}.  We first present an estimate which is valid for general polyhedral domains.
\begin{theorem}
\label{t3-3}
Let  $\Omega \subset \mathbb{R}^n$, $n=2,3$,  be an arbitrary polyhedral domain, and let $u_h \in S_0$ be a standard Lagrange finite element approximation defined on an arbitrary shape-regular simplicial decomposition of $\Omega$ having minimum mesh diameter $\hmin$.  Then for $0 < t_0 \le T$, 
\begin{equation}
\begin{split}
\|(u-u_h)(t_0)\|_{L_\infty(\Omega)} \le& \|u_0-u_h^0\|_{L_\infty(\Omega)}
\\&+C(\Omega) (\ln \hmin)^2 [\E_{\infty,0}(u_h(0), g(0))+\E_{\infty,0}(u_h(t_0),g(t_0))
 \\ &+ \| \E_{\infty,0}(u_{h,t},g_t) \|_{L_1((0,T))} ],
\end{split}
\nonumber
\end{equation}
where $E_{\infty,0}$ is the $L_\infty$-type residual estimator defined in (\ref{eq2-4-4}).  
\end{theorem}

{\it Proof:}
We proceed by bounding the terms in (\ref{eq3-2-1}) using Lemma \ref{lem2-4-1}.  Recalling (\ref{eq3-1-2}) and (\ref{eq3-1-3}), we find that 
\begin{equation}
\|(\R u_h-u_h)(t_0)\|_{L_\infty(\Omega)} \le C(\Omega) (\ln \hmin)^2 \E_{\infty,0}(u_h(t),g(t))
\nonumber
\end{equation}
and similarly
\begin{equation}
\|\R u_h(0)-u^0\|_{L_\infty(\Omega)} \le \|u_0-u_h^0\|_{L_\infty(\Omega)}+C(\Omega) (\ln \hmin)^2 \E_{\infty,0}(u_h^0,g(0)).
\nonumber
\end{equation}
Finally, 
\begin{equation}
\|(u_h-\R u_h)_t\|_{L_1([0,t], L_\infty(\Omega))} \le C(\Omega) (\ln \hmin)^2 \|\E_{\infty,0}(u_{h,t}, g_t)\|_{L_1([0,t])}.
\nonumber
\end{equation}
Inserting the above inequalities into (\ref{eq3-2-2}) completes the proof.  
\hfill $\Box$

We next present a theorem which allows us to bound the main error term with a higher-order estimator.  However, this estimate only holds for quadratic and higher-order elements and convex polygonal domains.  
\begin{theorem}
\label{t3-4}
Assume that $\Omega$ is a convex polyhedral domain in $\mathbb{R}^n$, $n=2,3$, with maximum vertex (if $n=2$) or edge (if $n=3$) opening angle $\frac{\pi}{\alpha}$, $\alpha>1$.  Let also $2 < p,q \le\infty$ satisfy $\frac{n}{2p}+\frac{1}{q}<\frac{1}{2}$ and $\frac{2}{\alpha-1}<p$.  Finally, assume that  the polynomial degree of the finite element space $S_0$ is at least two.  Then with $\hmin=\min_{K \in \T_h} \diam(K)$, 
\begin{equation}
\begin{split}
\|(u-&u_h)(t_0)\|_{L_\infty(\Omega)} \le \|u_0-u_h^0\|_{L_\infty(\Omega)}
\\&+C(\Omega) (\ln \hmin)^2 [\E_{\infty,0}(u_h(0), g(0))+\E_{\infty,0}(u_h(t_0), g(t_0))]
\\& +C_{p,q}(t_0) C(p,\Omega) \|\E_{p,-1}(u_{h,t}, g_t)\|_{L_q([0,t_0])}. 
\end{split}
\nonumber
\end{equation}
\end{theorem}
{\it Proof:}  The proof is completely analogous to that of Theorem \ref{t3-3} above, the only difference being that we now employ Lemma \ref{lem2-4-2} in addition to Lemma \ref{lem2-4-1}.  \hfill $\Box$

\section{Analysis of the fully discrete scheme}

\subsection{Fully discrete finite element approximation}

As in \S \ref{sec1}, let $0=t_0 < t_1 <...<t_N=T$, $I_i=(t_{i-1}, t_i)$, and $\tau_i=|I_i|$.  For each $0 \le i \le N$, let $\T_i$ be a shape-regular simplicial decomposition of $\Omega$.  Let also $S^i$ be a space of continuous piecewise polynomials of degree $k$ on $\T_i$, $k \ge 1$, and $S_0^i=S^i \cap H_0^1(\Omega)$.  Defining $v^i(x)=v(t_i,x)$ for $v \in C(\Omega \times [0,T])$, we discretize the weak form of (\ref{eq1-1-1}) using the backward Euler method as follows.  Let $u_h^0\in S_0^0$ be some approximation to $u_0$.  $u_h^i \in S_0^i$, $1 \le i \le N$, is then defined via the recursion
\begin{equation}
\left ( \frac{u_h^i-u_h^{i-1}}{\tau_i}, \phi_i \right ) + (\nabla u_h^i,\nabla  \phi_i) = (f^i, \phi_i)  \hbox{ for all } \phi_i \in S_0^i.
\label{eq4-1-1}
\end{equation}
In order to obtain a discrete approximation to $u$ on the whole parabolic domain $\Omega \times [0,T]$, we interpolate the functions $u_h^i$ linearly between $t_{i-1}$ and $t_i$:  
\begin{equation}
u_h(x,t)=(1-\frac{t-t_{i-1}}{\tau_i}) u_h^{i-1}(x)+\frac{t-t_{i-1}}{\tau_i} u_h^i(x), ~t_{i-1} \le t \le t_i.
\nonumber
\end{equation}
Finally, we define $u_{h,t}^i=\frac{\partial}{\partial t} u_h|_{I_i}$, that is, 
\begin{equation}
u_{h,t}^i(x)=\frac{u_h^i(x)-u_h^{i-1}(x)}{\tau_i}, ~ i \ge 1.
\label{eq4-1-4}
\end{equation}

Next we define $L_2$ projections onto $S^i$ and $S_0^i$.  For $0 \le i \le N$, we define $P_h^i:L_2(\Omega) \rightarrow S^i$ and $P_{h,0}^i :L_2(\Omega) \rightarrow S_0^i$ by
\begin{eqnarray}
(P_h^i u, v_i)&=&(u,v_i) ~ \forall ~v_i \in S^i,
\nonumber
\\ (P_{h,0}^i u,v_i)& = & (u, v_i) ~ \forall~ v_i \in S_0^i.
\nonumber
\end{eqnarray}
The discrete Laplacian $-\Delta_h^i:H_0^1(\Omega) \rightarrow S_0^i + P_h^i f^i $ is then given by
\begin{equation}
(-\Delta_h^i u, v_i)=(\nabla u,\nabla v_i) \hbox{ for all } v_i \in S_0^i.
\label{eq2-1-3}
\end{equation}

The weak-form fully discrete scheme (\ref{eq4-1-1}) may easily be transformed into the pointwise equation
\begin{equation}
\frac{u_h^i-P_{h,0}^i u_h^{i-1}}{\tau_i}-\Delta_h^i u_h^i=P_h^i f^i.
\nonumber
\end{equation}
Referring to (\ref{eq4-1-4}), we thus find that
\begin{equation}
u_{h,t}^i-\Delta_h^i u_h^i=P_h^i f^i+\frac{P_{h,0}^i u_h^{i-1}-u_h^{i-1}}{\tau_i}, ~ i \ge 1.
\label{eq2-1-5}
\end{equation}

\subsection{Elliptic reconstruction}

We define the elliptic reconstruction by first defining it at the time nodes and then interpolating linearly between them.   Following (\ref{eq3-1-2}), for $0 \le i \le N$ we let
\begin{equation}
(\nabla \R^i u_h^i,\nabla v)= (g^i, v) ~\forall ~ v \in H_0^1(\Omega),
\label{eq4-1-5b}
\end{equation}
where
\begin{equation}
g^i= \left \{ \begin{array}{rl} & -\Delta_h^0 u_h^0+f^0-P_h^0 f^0,~ i=0,
\\& f^i-u_{h,t}^i, ~i \ge 1.
 \end{array} \right .
\label{eq4-1-5}
\end{equation}
We then obtain the time-continuous elliptic reconstruction
\begin{equation}
\R u_h= (1-\frac{t-t_{i-1}}{\tau_i}) \R^{i-1} u_h^{i-1}+\frac{t-t_{i-1}}{\tau_i} \R^i u_h^i, ~t_{i-1} \le t \le t_i.
\nonumber
\end{equation}

Using the formula (\ref{eq2-1-5}), it is easy to compute that 
\begin{equation}
f^i-u_{h,t}^i =[-\Delta_h^i u_h^i]+[f^i-P_h^i f^i]+[\frac{u_h^{i-1}-P_{h,0}^i u_h^{i-1}}{\tau_i}].  
\label{eq4-1-8}
\end{equation}
This relationship combined with the definitions (\ref{eq2-1-3}) and (\ref{eq4-1-5}) yields the Galerkin orthogonality relationship
\begin{equation}
(\nabla( u_h^i-\R^i u_h^i), \nabla v_h)=0 ~ \forall ~ v_h \in S_0^i.
\label{eq4-1-5a}
\end{equation}
Similarly, we have on $I_i$ that
\begin{equation}
(\nabla (u_{h,t}^i-(\R u_h)_t),\nabla v_h)=0 ~ \forall  ~v_h \in S_0^{i-1} \cap S_0^i.
\label{eq4-1-5c}
\end{equation}

We finally state an error equation which will play a fundamental role in our analysis.  For $t \in I_i$ and $ \phi \in H_0^1(\Omega)$, 
\begin{equation}
\begin{split}
((u-\R u_h)_t,\phi)+(\nabla(u-\R u_h),\nabla \phi)=& ((u_h-\R u_h)_t,\phi)+(f-f^i,\phi)
\\ &+(1-\frac{t-t_{i-1}}{\tau_i})(g^i-g^{i-1},\phi).
\end{split}
\label{eq4-1-7}
\end{equation}
Comparing (\ref{eq4-1-7}) with (\ref{eq3-1-3a}), we see that (\ref{eq4-1-7}) contains additional terms which all result from the time discretization of (\ref{eq1-1-1}).  

\begin{remark}{\rm
From (\ref{eq4-1-8}), we see that the elliptic reconstruction for the fully discrete problem lifts the sum of a {\it discrete Laplacian term}, a {\it spatial data approximation term}, and a {\it mesh coarsening term}.  Note that the last term is nonzero only if $S_0^{i-1} \nsubseteq S_0^{i}$, that is, if the mesh is coarsened in the $i$-th time step.  Our definitions of $\R^i$ for $i=0$ and $ i > 0$ therefore differ only by the exclusion of the mesh coarsening term when $i=0$.   The  definition of $\R^i$ in \cite{LM06} does not incorporate the data approximation and mesh coarsening terms in (\ref{eq4-1-8}) into the right hand side of the equation solved by $\R^i u_h^i$.  The advantage of including these terms is that the error equation (\ref{eq4-1-7}) now includes only the term $((u_h -\R u_h)_t, \phi)$ plus terms which result from the discretization in time, resulting finally in a posteriori estimators which have a simpler structure. 
 }\end{remark}


\subsection{Fully discrete reconstruction results}
Here we present three alternative results.  
\begin{theorem}
\label{t4-1}
Let all assumptions and definitions be as in \S 4.1 and \S 4.2.  Then for any $1 \le j \le N$, 
\begin{equation}
\begin{split}
\|(u-u_h)(t_j)\|_{L_\infty(\Omega)} \le& \|(\R u_h-u_h)(t_j)\|_{L_\infty(\Omega)}+\|(u-\R u_h)(0)\|_{L_\infty(\Omega)}
\\ & +\|(u_h -\R u_h)_t \|_{L_1([0,t_j]; L_\infty(\Omega))}
\\ &+ \sum_{i=1}^j \int_{I_i} \|f-f^{i}\|_{L_\infty(\Omega)} \d t+\frac{\tau_i}{2}\|g^i-g^{i-1}\|_{L_\infty(\Omega)}. 
\end{split}
\label{eq4-2-2}
\end{equation}
Alternatively,  for $2 < p,q \le \infty$ satisfying $\frac{n}{2p}+\frac{1}{q}<\frac{1}{2}$ and for any $i \le j \le N$, we have
\begin{equation}
\begin{split}
\|(u-u_h)(t_j)\|_{L_\infty(\Omega)} \le& \|(\R u_h-u_h)(t_j)\|_{L_\infty(\Omega)}+\|(u-\R u_h)(0)\|_{L_\infty(\Omega)}
\\ &+C_{p,q}(t_j) \|(u_h -\R u_h)_t\|_{L_q([0,t_j]; W_p^{-1}(\Omega))}
\\ &+  \sum_{i=1}^j \int_{I_i} \|f-f^{i}\|_{L_\infty(\Omega)} \d t+\frac{\tau_i}{2}\|g^i-g^{i-1}\|_{L_\infty(\Omega)}. 
\end{split}
\label{eq4-2-1}
\end{equation}
Finally let $c(n)$ be the constant from $(\ref{eq2-3-5})$.  Then
\begin{equation}
\begin{split}
\|(u-&u_h)(t_j)\|_{L_\infty(\Omega)}  \le  \|u_0-u_{h0}\|_{L_\infty(\Omega)} 
\\ &+ (2+c(n) \ln \frac{t_j}{\tau_j}) \max_{0 \le i \le j} \|\R^i u_h^i -u_h^i\|_{L_\infty(\Omega)}
\\ &+  \sum_{i=1}^j \int_{I_i} \|f-f^{i}\|_{L_\infty(\Omega)} \d t+\frac{\tau_i}{2}\|g^i-g^{i-1}\|_{L_\infty(\Omega)}. 
\end{split}
\label{eq4-2-2a}
\end{equation}
\end{theorem}

{\it Proof:}  
We first fix a point  $x_0 \in \Omega$ with $\|(u-u_h)(t_j)\|_{L_\infty(\Omega)}=|(u-u_h)(x_0, t_j)|$ and compute 
\begin{equation}
|(u-u_h)(x_0,t_j)| \le |(u-\R u_h )(x_0,t_j)|+\|(\R u_h -u_h)(t_j)\|_{L_\infty(\Omega)}.  
\nonumber
\end{equation}
Using (\ref{eq2-3-1}), (\ref{eq4-1-7}), and (\ref{eq2-3-4}), we find that
\begin{equation}
\begin{split}
(u-&\R u_h)(x_0,t_j)=\int_\Omega G(x_0,t_j; y,0) (u-\R u_h)(y,0) \d y
\\ & +\int_0^{t_j} \int_\Omega G(x_0, t_j; y,s) (u_h -\R u_h)_s (y,s) \d y \d s
\\ & + \sum_{i=1}^j \int_{I_i} \int_\Omega G(x_0, t_j; y,s) (f-f^i)(y,s) \d y \d s
\\&+ \sum_{i=1}^j \int_{I_i} \int_\Omega G(x_0, t_j; y,s) (1-\frac{s-t_{i-1}}{\tau_i}) (g^i-g^{i-1}) \d y \d s
\\  \le &  \|(u-\R u_h)(0)\|_{L_\infty(\Omega)} + \int_0^{t_j} \int_\Omega G(x_0, t_j; y,s) (u_h -\R u_h)_s (y,s) \d y \d s
\\ &+ \sum_{i=1}^j \int_{I_i} \|f-f^{i}\|_{L_\infty(\Omega)} \d t+\frac{\tau_i}{2}\|g^i-g^{i-1}\|_{L_\infty(\Omega)}. 
\end{split}
\label{eq4-2-4}
\end{equation}
The term $\int_0^{t_j} \int_\Omega G(x_0, t_j; y,s) (u_h -\R u_h)_t (y,s) \d y \d s$ above may then be bounded by using (\ref{eq2-3-4}) in order to obtain (\ref{eq4-2-2}), or by using (\ref{eq2-3-3}) to obtain (\ref{eq4-2-1}).  

In order to prove (\ref{eq4-2-2a}), we begin as in (\ref{eq4-2-4}), split the second integral into two integrals, and perform integration by parts in time on $[0, t_{j-1}]$ to compute
\begin{equation}
\begin{split}
(u-&\R u_h)(x_0,t_j)=[\int_\Omega G(x_0,t_j; y,0) (u-\R u_h)(y,0) \d y]
\\ & +[\int_0^{t_j} \int_\Omega G(x_0, t_j; y,s) (u_h -\R u_h)_s (y,s) \d y \d s ]
\\ & + [ \sum_{i=1}^j \int_{I_i} \int_\Omega G(x_0, t_j; y,s) (f-f^i)(y,s) \d y \d s ]
\\&+ [ \sum_{i=1}^j \int_{I_i} \int_\Omega G(x_0, t_j; y,s) (1-\frac{s-t_{i-1}}{\tau_i}) (g^i-g^{i-1}) \d y \d s ]
\\  =&  [I]+[II+III+IV+V]+[VI]+[VII],
\end{split}
\label{eq4-2-5}
\end{equation}
where
\begin{eqnarray*}
I&=&\int_\Omega G(x_0,t_j; y,0) (u-\R u_h)(y,0) \d y,
\\II&=&  \int_{\Omega} G(x_0, t_j; y,t_{j-1}) (u_h -\R u_h)(y, t_{j-1}) \d y,
\\ III&=&  -  \int_{\Omega} G(x_0, t_j; y,0) (u_h -\R u_h)(y, 0) \d y,
 \\ IV&=& - \int_0^{t_{j-1}}  \int_\Omega G_s(x_0, t_j; y,s)  (u_h -\R u_h) (y,s) \d y \d s,
\\ V&=&  \int_{t_{j-1}}^{t_j}  \int_\Omega G(x_0, t_j; y,s) (u_h -\R u_h)_s (y,s) \d y \d s,
\\ VI &=&  \sum_{i=1}^j \int_{I_i} \int_\Omega G(x_0, t_j; y,s) (f-f^i)(y,s) \d y \d s,
\\ VII &=&  \sum_{i=1}^j \int_{I_i} \int_\Omega G(x_0, t_j; y,s) (1-\frac{s-t_{i-1}}{\tau_i}) (g^i-g^{i-1}) \d y \d s.
\end{eqnarray*}

In order to bound the terms in (\ref{eq4-2-5}), we use (\ref{eq2-3-3}) to find that
\begin{equation}
\begin{split}
I +II+III=& \int_\Omega G(x_0, t_j; y,0) (u-u_h)(y,0) \d y
\\ &+ \int_\Omega G(x_0, t_j; y,t_{j-1}) (u_h -\R u_h)(y, t_{j-1}) \d y
\\  \le & \|G(x_0, t_j; \cdot, 0)\|_{L_1(\Omega)} \|u_0 -u_{h0}\|_{L_\infty(\Omega)} 
\\ & + \|G(x_0, t_j; \cdot, t_{j-1})\|_{L_1(\Omega)} \|(u_h -\R u_h)(t_{j-1}) \|_{L_\infty(\Omega)}
\\  \le & \|u_0 -u_{h0}\|_{L_\infty(\Omega)} + \|(u_h -\R u_h)(t_{j-1}) \|_{L_\infty(\Omega)}.
\end{split}
\label{eq4-2-6}
\end{equation}
Employing (\ref{eq2-3-5}), we obtain
\begin{equation}
\begin{split}
IV  \le &  \int_{0}^{t_{j-1}} \|G_s(x_0, t_0; \cdot, s)\|_{L_1(\Omega)} \|(u_h -\R u_h)(s)\|_{L_\infty(\Omega)} \d s
\\  \le & \|u_h -\R u_h \|_{L_\infty(\Omega \times (0, t_{j-1}))} \int_0^{t_{j-1}} \frac{c(n)}{t_j-s} \d s
\\  \le & c(n) \ln \frac{ t_j}{\tau_j} \max_{0 \le i \le j-1} \| u_h^i-\R^i u_h^i\|_{L_\infty(\Omega)} .
\end{split}
\label{eq4-2-7}
\end{equation}
In order to bound the term $V$, we calculate
\begin{equation}
\begin{split}
V  \le & \|G(x_0, t_0, \cdot, \cdot)\|_{L_1(\Omega \times I_j)} \|(u_h - \R u_h)_t\|_{L_\infty(\Omega \times I_j)}
\\  \le & \tau_j \|\frac{(u_h^j-\R^j u_h^j)-(u_h^{j-1}-\R^{j-1} u_h^{j-1})}{\tau_j}\|_{L_\infty(\Omega)}
\\  \le & \|u_h^j-\R^j u_h^j\|_{L_\infty(\Omega)}+ \|u_h^{j-1}-\R^{j-1} u_h^{j-1}\|_{L_\infty(\Omega)}.
\end{split}
\nonumber
\end{equation}
Finally, we compute directly that
\begin{equation}
VI+VII \le  \sum_{i=1}^j \int_{I_i} \|f-f^{i}\|_{L_\infty(\Omega)} \d t+\frac{\tau_i}{2}\|g^i-g^{i-1}\|_{L_\infty(\Omega)}.
\label{eq4-2-9}
\end{equation}
Collecting the previous inequalities and inserting them into (\ref{eq4-2-5}) completes the proof of (\ref{eq4-2-2a}).  
\hfill $\Box$

\begin{remark}{\rm One may approach the proofs of (\ref{eq4-2-2}) and (\ref{eq4-2-2a}) of Theorem \ref{t4-1} and (\ref{eq3-2-2}) of Theorem \ref{t3-1} from the perspective of semigroup theory instead of using fundamental solutions.  Let $E$ be the semigroup generated by the Laplace operator, i.e., let $E(t)v_0$ be the solution of $v_t-\Delta v=0, v(0)=v_0$.  Then using Duhamel's principle we have
\begin{equation}
u(t)=E(t) u_0+\int_0^t E(t-s) f(s) \d s.
\nonumber
\end{equation}
Assume that the stability and strong stability properties 
\begin{eqnarray}
\|E(t)\|_{L_\infty \rightarrow L_\infty} &\le& 1,~t>0,
\label{eq4-2-11}
\\ \|E'(t)\|_{L_\infty \rightarrow L_\infty} & \le& \frac{C}{t}, ~t>0
\label{eq4-2-12}
\end{eqnarray}
hold.  The bounds (\ref{eq3-2-2}) and (\ref{eq4-2-2a}) are then easily obtained by respectively combining (\ref{eq3-1-3a}) and (\ref{eq4-1-7}) with (\ref{eq4-2-11}).  (\ref{eq4-2-2a}) may be obtained by combining (\ref{eq4-2-11}), (\ref{eq4-2-12}), and (\ref{eq4-1-7}).  We note, however, that we are not aware of a proof for the results (\ref{eq4-2-11}) and (\ref{eq4-2-12}) under the weak restrictions that we have placed on $\Omega$.  The standard reference \cite{St74}, for example, assumes that $\partial \Omega$ is $C^2$ in order to obtain the analyticity of $E$ in $C^0$ and thus obtain (\ref{eq4-2-12}).  
}
\end{remark}

We finally note that it is possible to sharpen the reconstruction estimate (\ref{eq4-2-2a}) somewhat, though at the cost of a more complex result.  In particular, we may employ the estimates (\ref{eq2-3-8}) and (\ref{eq2-3-9}) that reflect the dissipation of the heat kernel when $n=2$ and also accumulate the spatial errors in $\ell_1$ instead of $\ell_\infty$.  Let $\phi_i(s)$ be the piecewise linear ``hat'' function satisfying $\phi_i(t_i)=1$, $\phi_i(t_m)=0$ for $m \neq i$.   Instead of (\ref{eq4-2-7}), we may then calculate 
\begin{equation}
\begin{split}
IV \le & \int_{0}^{t_{j-1}}   \|G_s(x_0, t_0; \cdot, s)\|_{L_1(\Omega)} \|(u_h -\R u_h)(s)\|_{L_\infty(\Omega)} \d s
\\ \le & c(n) [\|u_h^0-\R^0 u_h^0 \|_{L_\infty(\Omega)} \int_0^{t_1} \phi_0(s) \frac{1}{t_j-s} (1-e^{-\frac{\diam(\Omega)^2}{9(t_j-s)}}) \d s
\\ & +\sum_{i=1}^{j-2} \|u_h^i-\R^i u_h^i\|_{L_\infty(\Omega)} \int_{t_{i-1}}^{t_{i+1}} \phi_i(s) \frac{1}{t_j-s} (1-e^{-\frac{\diam(\Omega)^2}{9(t_j-s)}}) \d s 
\\ & +\|u_h^{j-1}-\R^{j-1} u_h^{j-1}\|_{L_\infty(\Omega)} \int_{t_{j-2}}^{t_{j-1}} \phi_{j-1}(s) \frac{1}{t_j-s} (1-e^{-\frac{\diam(\Omega)^2}{9(t_j-s)}}) \d s ]
\\ \le & c(n) [ \frac{\tau_1}{2(t_j-t_1)} (1-e^{-\frac{\diam(\Omega)^2}{9(t_j-t_1)}}) \|u_h^0-\R^0 u_h^0 \|_{L_\infty(\Omega)} 
\\ & +\sum_{i=1}^{j-2} \frac{1}{2}(\frac{\tau_i}{t_j-t_i}+\frac{\tau_{i+1}}{t_j-t_{i+1}}) (1-e^{-\frac{\diam(\Omega)^2}{9(t_j-t_{i+1})}})  \|u_h^i-\R^i u_h^i\|_{L_\infty(\Omega)}
\\ &+ \frac{\tau_{j-1}}{2 \tau_j}(1-e^{-\frac{\diam(\Omega)^2}{9 \tau_j}}) \|u_h^{j-1}-\R^{j-1} u_h^{j-1}\|_{L_\infty(\Omega)}].  
\end{split}
\label{eq4-3-13}
\end{equation}
Employing (\ref{eq4-3-13}) instead of (\ref{eq4-2-7}) and similarly inserting (\ref{eq2-3-8}) into (\ref{eq4-2-6}) and (\ref{eq4-2-9})  leads to the following result.
\begin{prop}
\label{cor4-4}
If the spatial dimension $n=2$ and the conditions of Theorem \ref{t4-1} are met, then
\begin{equation}
\begin{split}
\|& (u-u_h)(t_j)\|_{L_\infty(\Omega)} \le  \beta_4 (0) \|u_0-u_{h0}\|_{L_\infty(\Omega)}
\\ &+c(n)   \frac{\tau_1}{2(t_j-t_1)} \beta_9 (t_1) \|u_h^0-\R^0 u_h^0 \|_{L_\infty(\Omega)}  
\\ & c(n) \sum_{i=1}^{j-2} \frac{1}{2}(\frac{\tau_i}{t_j-t_i}+\frac{\tau_{i+1}}{t_j-t_{i+1}}) \beta_9 (t_{i+1}) \|u_h^i-\R^i u_h^i\|_{L_\infty(\Omega)}
\\ & +(2+c(n) \frac{\tau_{j-1}}{2 \tau_j}  \beta_9 (t_{j-1}) ) \|u_h^{j-1}-\R^{j-1} u_h^{j-1}\|_{L_\infty(\Omega)}
 +2 \|u_h^j-\R^j u_h^j\|_{L_\infty(\Omega)} 
\\ & + \sum_{i=1}^j \int_{I_i} \beta_4 (t) \|f-f^{i}\|_{L_\infty(\Omega)} \d t
+\sum_{i=1}^j \frac{\tau_i}{2}\beta_4 (t_i) \|g^i-g^{i-1}\|_{L_\infty(\Omega)}. 
\label{eq4-3-14}
\end{split}
\end{equation}
Here
\begin{eqnarray}
\beta_4(s)&=&1-e^{-\frac{\diam(\Omega)^2}{4(t_j-s)}},
\nonumber
\\ \beta_9(s)&=&1-e^{-\frac{\diam(\Omega)^2}{9(t_j-s)}}.
\nonumber
\end{eqnarray}
\end{prop} 

In contrast to (\ref{eq4-2-2a}), (\ref{eq4-3-14}) reflects damping of the error effects from times $t_i$ for which $t_j-t_i >>0$.  This damping is reflected in two ways.  First, in (\ref{eq4-3-14}) the spatial indicators $\|u_h^i-\R u_h^i\|_{L_\infty(\Omega)}$ accumulate in $\ell_1$ instead of in $\ell_\infty$ as in (\ref{eq4-2-2a}).  The weights  $\frac{1}{2}(\frac{\tau_i}{t_j-t_i}+\frac{\tau_{i+1}}{t_j-t_{i+1}})$ in (\ref{eq4-3-14}) deemphasize the spatial indicators for times $t_i << t_j$ relative to contributions from times $t_i \approx t_j$.  This yields a sharper bound than (\ref{eq4-2-2a}), where the spatial indicators are all weighted equally. Note that similar weights can also be obtained when using a duality-reconstruction combination (cf. \cite{LM07}).   In addition, the dissipation estimates (\ref{eq2-3-8}) and (\ref{eq2-3-9}) are used above to quantify damping in the initial data error (first line above) and time discretization indicators (last line) as well as in the spatial errors.  While not dramatic, these dissipation effects can be substantial.  If for example $t_j=\diam(\Omega)=1$, then $\beta_4(0) \approx .222$, corresponding roughly to one level of refinement in the initial mesh if piecewise linear elements are used.  

We finally note that (\ref{eq4-3-14}) still holds when $n>2$, but the dissipation weights $\beta_4(s)$ and $\beta_9(s)$ take a different form.

\subsection{A posteriori error estimates for the fully discrete problem}
We finally obtain three different a posteriori estimates.

\begin{theorem}
\label{t4-2}
Let $\Omega \subset \mathbb{R}^n$, $n=2,3$, be an arbitrary polyhedral domain, and define $\hmin=\min_{0 \le i \le N} \min_{T \in \T_i} h_T$.  Then for any $1 \le j \le N$, 
\begin{equation}
\begin{split}
\|(u-u_h)(t_j)\|_{L_\infty(\Omega)} \le & \|u_0-u_h^0\|_{L_\infty(\Omega)}
\\ &+ C(\Omega) (\ln \hmin)^2[ \E_{\infty,0} (u_h^j, g_j)+\E_{\infty,0}(u_h^0, g_0)
\\ & +\sum_{i=1}^j \hat{\E}_{\infty,0} (u_h^i-u_h^{i-1}, g^i-g^{i-1}; \T_{i-1}, \T_i)]
\\ &+ \sum_{i=1}^j \int_{I_i} \|f-f^{i}\|_{L_\infty(\Omega)} \d t+\frac{\tau_i}{2}\|g^i-g^{i-1}\|_{L_\infty(\Omega)}. 
\end{split}
\label{eq4-4-1}
\end{equation}
\end{theorem}

\begin{theorem}
\label{t4-3}
Assume that $\Omega$ is a convex polyhedral domain in $\mathbb{R}^n$, $n=2,3$ with maximum vertex (if $n=2$) or edge (if $n=3$) opening angle $\frac{\pi}{\alpha}$, $\alpha>1$.  Let also $2 < p,q \le\infty$ satisfy $\frac{n}{2p}+\frac{1}{q}<\frac{1}{2}$ and $\frac{2}{\alpha-1}<p$.  Finally, assume that  the polynomial degree of the finite element space $S_0$ is at least two.  Then for $1 \le j \le N$,
\begin{equation}
\begin{split}
\|(u-&u_h)(t_j)\|_{L_\infty(\Omega)} \le \|u_0-u_h^0\|_{L_\infty(\Omega)}
\\ & + C(\Omega) (\ln \hmin)^2 [\E_{\infty,0} (u_h^j, f^j-u_{h,t}^j)+\E_{\infty,0}(u_h^0, g_0)]
\\ & + C_{p,q}(t_j) C_p(\Omega) \left ( \sum_{i=1}^j \tau_i^{-q+1} \hat{\E}_{p,1} (u_h^i-u_h^{i-1}, g^i-g^{i-1}; \T_{i-1}, \T_i)^q\right )^{1/q}
\\ &+ \sum_{i=1}^j \int_{I_i} \|f-f^{i}\|_{L_\infty(\Omega)} \d t+\frac{\tau_i}{2}\|g^i-g^{i-1}\|_{L_\infty(\Omega)}. 
\end{split}
\label{eq4-4-2}
\end{equation}
\end{theorem}

\begin{theorem}
\label{t4-4}
Let $\Omega \subset \mathbb{R}^n$, $n=2,3$, be an arbitrary polyhedral domain, and define $\hmin=\min_{0 \le i \le N} \min_{T \in \T_i} h_T$.  Then for any $1 \le j \le N$,
\begin{equation}
\begin{split}
\|(u-u_h)(t_j)\|_{L_\infty(\Omega)} \le & \|u_0-u_h^0\|_{L_\infty(\Omega)}
\\ & +(2+c(n) \ln \frac{t_j}{\tau_j}) C(\Omega) (\ln \hmin)^2 \max_{0 \le i \le j} \E_{\infty,0}(u_h^i, g^i)
\\ &+ \sum_{i=1}^j \int_{I_i} \|f-f^{i}\|_{L_\infty(\Omega)} \d t+\frac{\tau_i}{2}\|g^i-g^{i-1}\|_{L_\infty(\Omega)},
\end{split}
\label{eq4-4-3}
\end{equation}
where $c(n)=\frac{3^n}{2^{n/2+1}}$.  
\end{theorem}

{\it Proof:}  The proofs of the three preceding theorems follow easily by inserting the estimates of Lemma \ref{lem2-4-1}, Lemma \ref{lem2-4-2}, Corollary \ref{cor2-7}, and Corollary \ref{cor2-8} into the estimates (\ref{eq4-2-2}), (\ref{eq4-2-1}), and (\ref{eq4-2-2a}) while recalling the definitions (\ref{eq4-1-5b}) and (\ref{eq4-1-5}) and the Galerkin orthogonality results (\ref{eq4-1-5a}) and (\ref{eq4-1-5c}).  \hfill $\Box$

\begin{remark}{\rm 
In contrast to (\ref{eq4-4-1}) and (\ref{eq4-4-2}), the estimate (\ref{eq4-4-3}) does not require the computation of residual-based estimators with respect to a finest common coarsening.   The fact that spatial errors in (\ref{eq4-4-3}) accumulate in $L_\infty$ in time as well as in space also is practically advantageous as it is easier to ensure that these errors are of the correct size at each time step in an adaptive algorithm.  In addition, (\ref{eq4-4-3}) more readily lends itself to establishing a convenient bound for $\|u-u_h\|_{L_\infty(\Omega \times [0,t_j])}$.  Thus for practical purposes, (\ref{eq4-4-3}) is of greatest interest among the results in this paper.  We also emphasize that (\ref{eq4-4-3}) does not include any unknown time-dependent constants. In fact, {\it all} unknown constants in (\ref{eq4-4-3}) stem from the use of a posteriori error estimators for {\it elliptic} problems.   

While not as practically advantageous as Theorem \ref{t4-4}, Theorem \ref{t4-2} and Theorem \ref{t4-3} also have interesting theoretical features.  (\ref{eq4-4-1}) is interesting in that it includes no time-dependent constants of any sort.  (\ref{eq4-4-2}) bounds the spatial errors at intermediate times in a weaker negative norm in which $u-u_h$ is of higher order for quadratic and higher-order finite element spaces. 
} \end{remark}

\begin{remark}{\rm One may sharpen Theorem \ref{t4-4} as in Corollary \ref{cor4-4}.  The latter result may easily be adapted to obtain an a posteriori estimate, but the resulting estimator is cumbersome and we thus do not record it here.}\end{remark}

\begin{thebibliography}{10}

\bibitem[Ar68]{Ar68}
{\sc D.~G. Aronson}, {\em Non-negative solutions of linear parabolic
  equations}, Ann. Scuola Norm. Sup. Pisa (3), 22 (1968), pp.~607--694.

\bibitem[Bo00]{Bo00}
{\sc M.~Boman}, {\em On a posteriori error analysis in the maximum norm}, PhD
  thesis, Chalmers University of Technology and G{\"o}teborg University, 2000.

\bibitem[Da00]{Dan00}
{\sc D.~Daners}, {\em Heat kernel estimates for operators with boundary
  conditions}, Math. Nachr., 217 (2000), pp.~13--41.

\bibitem[DDP00]{DDP00}
{\sc E.~Dari, R.~G. Dur{\'a}n, and C.~Padra}, {\em Maximum norm error
  estimators for three-dimensional elliptic problems}, SIAM J. Numer. Anal., 37
  (2000), pp.~683--700.

\bibitem[Da88]{Da88}
{\sc M.~Dauge}, {\em Elliptic boundary value problems on corner domains},
  vol.~1341 of Lecture Notes in Mathematics, Springer-Verlag, Berlin, 1988.

\bibitem[Da92]{Da92}
\leavevmode\vrule height 2pt depth -1.6pt width 23pt, {\em Neumann and mixed
  problems on curvilinear polyhedra}, Integral Equations Operator Theory, 15
  (1992), pp.~227--261.

\bibitem[Dav97]{Dav97}
{\sc E.~B. Davies}, {\em Non-{G}aussian aspects of heat kernel behaviour}, J.
  London Math. Soc. (2), 55 (1997), pp.~105--125.

\bibitem[DDE05]{DDE05}
{\sc K.~Deckelnick, G.~Dziuk, and C.~M. Elliott}, {\em Computation of geometric
  partial differential equations and mean curvature flow}, Acta Numer., 14
  (2005), pp.~139--232.

\bibitem[De06]{De06}
{\sc A.~Demlow}, {\em Localized pointwise a posteriori error estimates for
  gradients of piecewise linear finite element approximations to second-order
  quasilinear elliptic problems}, SIAM J. Numer. Anal., 44 (2006), pp.~494--514.

\bibitem[De07]{Dem07}
\leavevmode\vrule height 2pt depth -1.6pt width 23pt, {\em Local a posteriori
  estimates for pointwise gradient errors in finite element methods for
  elliptic problems}, Math. Comp., 76 (2007), pp.~19--42.

\bibitem[EJ95]{EJ95}
{\sc K.~Eriksson and C.~Johnson}, {\em Adaptive finite element methods for
  parabolic problems. {II}. {O}ptimal error estimates in {$L\sb \infty L\sb 2$}
  and {$L\sb \infty L\sb \infty$}}, SIAM J. Numer. Anal., 32 (1995),
  pp.~706--740.

\bibitem[LM06]{LM06}
{\sc O.~Lakkis and C.~Makridakis}, {\em Elliptic reconstruction and a
  posteriori error estimates for fully discrete linear parabolic problems},
  Math. Comp., 75 (2006), pp.~1627--1658.

\bibitem[LM07]{LM07}
\leavevmode\vrule height 2pt depth -1.6pt width 23pt, {\em A posteriori error
  estimates for parabolic problems via elliptic reconstruction and duality}.
\newblock arXiv, 2007.

\bibitem[LN03]{LN03}
{\sc X.~Liao and R.~H. Nochetto}, {\em Local a posteriori error estimates and
  adaptive control of pollution effects}, Numer. Methods Partial Differential
  Equations, 19 (2003), pp.~421--442.

\bibitem[MN03]{MN03}
{\sc C.~Makridakis and R.~H. Nochetto}, {\em Elliptic reconstruction and a
  posteriori error estimates for parabolic problems}, SIAM J. Numer. Anal., 41
  (2003), pp.~1585--1594.

\bibitem[Noc95]{Noc95}
{\sc R.~H. Nochetto}, {\em Pointwise a posteriori error estimates for elliptic
  problems on highly graded meshes}, Math. Comp., 64 (1995), pp.~1--22.

\bibitem[NSSV06]{NSSV06}
{\sc R.~H. Nochetto, A.~Schmidt, K.~G. Siebert, and A.~Veeser}, {\em Pointwise
  a posteriori error estimates for monotone semilinear problems}, Numer. Math.,
  104 (2006), pp.~515--538.

\bibitem[NSV03]{NSV03}
{\sc R.~H. Nochetto, K.~G. Siebert, and A.~Veeser}, {\em Pointwise a posteriori
  error control for elliptic obstacle problems}, Numer. Math., 95 (2003),
  pp.~163--195.

\bibitem[NSV05]{NSV05}
\leavevmode\vrule height 2pt depth -1.6pt width 23pt, {\em Fully localized a
  posteriori error estimators and barrier sets for contact problems}, SIAM J.
  Numer. Anal., 42 (2005), pp.~2118--2135.

\bibitem[SS05]{SS05}
{\sc A.~Schmidt and K.~G. Siebert}, {\em Design of adaptive finite element
  software}, vol.~42 of Lecture Notes in Computational Science and Engineering,
  Springer-Verlag, Berlin, 2005.
\newblock The finite element toolbox ALBERTA, With 1 CD-ROM (Unix/Linux).

\bibitem[St74]{St74}
{\sc H.~B. Stewart}, {\em Generation of analytic semigroups by strongly
  elliptic operators}, Trans. Amer. Math. Soc., 199 (1974), pp.~141--162.

\bibitem[Th97]{Th97}
{\sc V.~Thom{\'e}e}, {\em Galerkin finite element methods for parabolic
  problems}, vol.~25 of Springer Series in Computational Mathematics,
  Springer-Verlag, Berlin, 1997.

\bibitem[Wh73]{Wh73}
{\sc M.~F. Wheeler}, {\em {A priori L$_2$ error estimates for Galerkin
  approximations to parabolic partial differential equations.}}, SIAM J. Numer.
  Anal., 10 (1973), pp.~723--759.

\end{thebibliography}

\end{document}